# Osculating direction curves and their applications

Mehmet Önder[a], Sezai Kızıltuğ[b]

[a]Celal Bayar University, Faculty of Arts and Sciences, Department of Mathematics, Muradiye Campus, Muradiye, Manisa, Turkey.
E-mail: mehmet.onder@cbu.edu.tr

[b]Erzincan University, Faculty of Arts and Sciences, Department of Mathematics, , Erzincan, Turkey.
E-mail: skiziltug24@hotmail.com

**Abstract**
In this study, we define a new type of direction curves in the Euclidean 3-space such as osculating-direction curve. We give the characterizations for these curves. Moreover, we obtain the relationships between osculating direction curves and some special curves such as helix, slant helix or rectifying curves.

**MSC:** 53A04

**Key words:** Associated curve; osculating-direction curve, osculating-donor curve.

## 1. Introduction

The most fascinating and important subject of curve theory is to obtain the characterizations for a regular curve or a curve pair. Helix, slant helix, plane curve, spherical curve, etc. are the well-known examples of special curves [1,9,10,17,20] and these special curves, especially the helices, are used in many applications [2,7,8,16]. Furthermore, other definitions of special curves can be given by considering Frenet planes. The curve for which the position vector always lies on its rectifying, osculating or normal planes are called rectifying curve, osculating curve or normal curve, respectively [4]. In the Euclidean space $E^3$, rectifying, normal and osculating curves satisfy Cesaro's fixed point condition [15]. Namely, rectifying, normal and osculating planes of such curves always contain a particular point. It is well known that if all the normal or osculating planes of a curve in $E^3$ pass through a particular point, then the curve lies in a sphere or is a plane curve, respectively. Furthermore, if all rectifying planes of a non-planar curve in $E^3$ pass through a particular point, then the ratio of torsion and curvature of such curve is a non-constant linear function. In particular, there exists a simple relationship between rectifying curves and Darboux vectors (centrodes), which play some important roles in mechanics, kinematics as well as in differential geometry in defining the curves of constant precession [5].

Moreover, the curve pairs for which there exits some relationships between their Frenet vectors or curvatures are another examples of special space curves. Involute-evolute curves, Bertrand curves, Mannheim curves are the well-known examples of curve pairs and studied by many mathematicians [3,11-13,18,19].

Recently, a new curve pair in the Euclidean 3-space $E^3$ has been defined by Choi and Kim [6]. They have considered a unit vector field $X$ defined in the Frenet basis of a Frenet curve $\alpha$ and an integral curve $\gamma$ of $X$ and given the definition and characterizations of principal-directional curve and principal-donor curve in $E^3$. They have also give some applications of these curves to helix, slant helix or plane curve.

In the present paper, we give a new definition of curve pair such as osculating-direction curve and osculating donor curve in $E^3$. We obtain some characterizations for these curves and show that osculating-direction curve and its osculating donor curve form a Mannheim pair. Moreover, we give some applications of osculating-direction curve to some special curves such as helix, slant helix, rectifying curve in $E^3$.

**2. Preliminaries**

We now recall some basic concepts on classical differential geometry of space curves and the definitions of general helix, slant helix in the Euclidean 3-space $E^3$.

A curve $\alpha : I \to \mathbb{R}^3$ with unit speed, is a general helix if there is some constant vector $u$, so that $\langle T, u \rangle = \cos \theta$ is constant along the curve, where $\theta \neq \pi/2$ and $T(s) = \alpha'(s)$ is a unit tangent vector of $\alpha$ at $s$. We define the curvature of $\alpha$ by $\kappa(s) = \|\alpha''(s)\|$. If $\kappa(s) \neq 0$, then the curve $\alpha$ is called Frenet curve and the unit principal normal vector $N(s)$ of the curve $\alpha$ at $s$ is given by $\alpha''(s) = \kappa(s) N(s)$. The unit vector $B(s) = T(s) \times N(s)$ is called the unit binormal vector of $\alpha$ at $s$. Then $\{T, N, B\}$ is called the Frenet frame of $\alpha$. For the derivatives of the Frenet frame, the following Frenet-Serret formulae hold:

$$\begin{bmatrix} T' \\ N' \\ B' \end{bmatrix} = \begin{bmatrix} 0 & \kappa & 0 \\ -\kappa & 0 & \tau \\ 0 & -\tau & 0 \end{bmatrix} \begin{bmatrix} T \\ N \\ B \end{bmatrix} \tag{1}$$

where $\tau(s)$ is the torsion of the curve $\alpha$ at $s$. It is well-known that the curve $\alpha$ is a general helix if and only if $\dfrac{\tau}{\kappa}(s) = \text{constant}$. If both $\kappa(s) \neq 0$ and $\tau(s)$ are constant, we call $\alpha$ as a circular helix. A curve $\alpha$ with $\kappa(s) \neq 0$ is called a slant helix if the principal normal lines of $\alpha$ make a constant angle with a fixed direction. Also, a slant helix $\alpha$ in $E^3$ is characterized by the differential equation of its curvature $\kappa$ and its torsion $\tau$ given by

$$\frac{\kappa^2}{\left(\kappa^2 + \tau^2\right)^{3/2}} \left(\frac{\tau}{\kappa}\right)' = \text{constant}.$$

(See [10]).

Now, we define some associated curves of a curve $\alpha$ in $E^3$ defined on an open interval $I \subset \mathbb{R}$. For a Frenet curve $\alpha : I \to \mathbb{R}^3$, consider a vector field $X$ given by

$$X(s) = u(s) T(s) + v(s) N(s) + w(s) B(s), \tag{2}$$

where $u$, $v$ and $w$ are arbitrary differentiable functions of $s$ which is the arc length parameter of $\alpha$. Let

$$u^2(s) + v^2(s) + w^2(s) = 1, \tag{3}$$

holds. In [6], Choi and Kim introduced the definition of $X$-direction curve and $X$-donor curve in $E^3$ as follows.

**Definition 2.1.** ([6]) Let $\alpha$ be a Frenet curve in Euclidean 3-space $E^3$ and $X$ be a unit vector field satisfying the equations (2) and (3). The integral curve $\beta : I \to \mathbb{R}^3$ of $X$ is called an $X$-direction curve of $\alpha$. The curve $\alpha$ whose $X$-direction curve is $\beta$ is called the $X$-donor curve of $\beta$ in $E^3$.

**Definition 2.2.** ([6]) An integral curve of principal normal vector $N(s)$ (resp. binormal vector $B(s)$) of $\alpha$ in (2) is called the principal-direction curve (resp. binormal-direction curve) of $\alpha$ in $E^3$.

**Remark 2.1.** ([6]) A principal-direction (resp. the binormal-direction) curve is an integral curve of $X(s)$ with $u(s) = w(s) = 0$, $v(s) = 1$ (resp. $u(s) = v(s) = 0$, $w(s) = 1$) for all $s$ in (2).

### 3. Osculating direction curve and osculating donor curve in $E^3$

In this section, we will give definition of osculating-direction curve and osculating donor curve in $E^3$. Besides, we obtain some theorems and results characterizing these curves.

**Definition 3.1.** Let $\alpha$ be a Frenet curve in $E^3$ and $X$ be a unit vector field lying on the osculating plane of $\alpha$ and defined by
$$X(s) = u(s)T(s) + v(s)N(s), \ u(s) \neq 0, \ v(s) \neq 0, \ X'(s) \notin Sp\{T, N\}. \tag{4}$$
The integral curve $\gamma: I \to \mathbb{R}^3$ of $X(s)$ is called an osculating-direction curve of $\alpha$. The curve $\alpha$ whose osculating-direction curve is $\gamma$ is called the osculating-donor curve in $E^3$.

From Definition 3.1 it is clear that the integral curve $\gamma: I \to \mathbb{R}^3$ of $X(s)$ is a unit speed curve. Then, without loss of generality we can take $s$ as the arc length parameter of $\gamma$.

***Theorem 3.1.*** Let $\alpha: I \to \mathbb{R}^3$ be a Frenet curve and an integral curve of $X(s) = u(s)T(s) + v(s)N(s)$ be the curve $\gamma: I \to \mathbb{R}^3$. Then, $\gamma$ is an osculating-direction curve of $\alpha$ if and only if
$$u(s) = \sin\left(\int \kappa ds\right) \neq 0, \quad v(s) = \cos\left(\int \kappa ds\right) \neq 0. \tag{5}$$
**Proof:** Since $\gamma$ is the osculating-direction curve of $\alpha$, from Definition (3.1) we have
$$X(s) = u(s)T(s) + v(s)N(s), \tag{6}$$
and
$$u^2(s) + v^2(s) = 1. \tag{7}$$
Differentiating (6) with respect to $s$ and by using the Frenet formulas, it follows
$$X'(s) = (u' - v\kappa)T + (v' + u\kappa)N + v\tau B. \tag{8}$$
Since $X'(s) \notin Sp\{T, N\}$, we have that $X'$ and $B$ are linearly dependent. Then from (8) we can write
$$\begin{cases} u' - v\kappa = 0, \\ v' + u\kappa = 0, \\ v\tau \neq 0. \end{cases} \tag{9}$$
The solutions of first two differential equations are
$$u(s) = \sin\left(\int \kappa ds\right) \neq 0, \quad v(s) = \cos\left(\int \kappa ds\right) \neq 0,$$
respectively, which completes the proof.

***Theorem 3.2.*** Let $\alpha: I \to \mathbb{R}^3$ be a Frenet curve. If $\gamma$ is the osculating-direction curve of $\alpha$, then $\gamma$ is a Mannheim curve of $\alpha$.

**Proof:** Since $\gamma$ is an integral curve of $X$, we have $\gamma' = X$. Denote the Frenet frame of $\gamma$ by $\{\bar{T}, \bar{N}, \bar{B}\}$. Differentiating $\gamma' = X$ with respect to $s$ and by using Frenet formulas we get
$$X' = \bar{T}' = \bar{\kappa}\bar{N}. \tag{10}$$
Moreover, we know that $X'$ and $B$ are linearly dependent. Then from (10) we get $\bar{N}$ and $B$ are linearly dependent, i.e, $\gamma$ is a Mannheim curve of $\alpha$.

**Theorem 3.3.** Let $\alpha : I \to \mathbb{R}^3$ be a Frenet curve. If $\gamma$ is the osculating direction curve of $\alpha$, then the curvature $\bar{\kappa}$ and the torsion $\bar{\tau}$ of the curve $\gamma$ are given as follows,

$$\bar{\kappa} = \tau \cos\left(\int \kappa ds\right), \quad \bar{\tau} = \tau \sin\left(\int \kappa ds\right). \tag{11}$$

**Proof:** From (8) and (10), we have

$$\bar{\kappa} \bar{N} = v \tau B. \tag{12}$$

By considering Eq. (12) and Eq. (5) we obtain

$$\bar{\kappa} \bar{N} = \tau \cos\left(\int \kappa ds\right) B, \tag{13}$$

which gives us

$$\bar{\kappa} = \tau \cos\left(\int \kappa ds\right). \tag{14}$$

Besides, from (13) and (14) we can write

$$\bar{N} = B. \tag{15}$$

Then we have

$$\bar{B} = \bar{T} \times \bar{N} = \cos\left(\int \kappa ds\right) T - \sin\left(\int \kappa ds\right) N. \tag{16}$$

Differentiating (16) with respect to $s$ we have

$$\bar{B}' = \tau \sin\left(\int \kappa ds\right) B. \tag{17}$$

Since $\bar{\tau} = -\langle \bar{B}', \bar{N} \rangle = -\langle \bar{B}', B \rangle$, Eq. (17) gives us

$$\bar{\tau} = \tau \sin\left(\int \kappa ds\right). \tag{18}$$

Hence, the proof is completed.

**Corollary 3.1.** Let $\gamma$ be the osculating-direction curve of the curve $\alpha$. Then

$$X = \bar{T} = \frac{\bar{\tau}}{\tau} T + \frac{\bar{\kappa}}{\tau} N, \quad \bar{N} = B, \quad \bar{B} = \frac{\bar{\kappa}}{\tau} T - \frac{\bar{\tau}}{\tau} N.$$

**Proof:** The proof is clear from Theorem 3.3.

**Theorem 3.4.** Let $\gamma$ be an osculating-direction curve of $\alpha$ with the curvature $\bar{\kappa}$ and the torsion $\bar{\tau}$. Then the curvature $\kappa$ and the torsion $\tau$ of $\alpha$ are given by

$$\tau = \sqrt{\bar{\kappa}^2 + \bar{\tau}^2}, \quad \kappa = \frac{\bar{\kappa}^2}{\bar{\kappa}^2 + \bar{\tau}^2}\left(\frac{\bar{\tau}}{\bar{\kappa}}\right)'.$$

**Proof:** From (14) and (18), we easily get

$$\tau = \sqrt{\bar{\kappa}^2 + \bar{\tau}^2}. \tag{19}$$

Substituting (19) into (14) and (18), we obtain

$$\cos\left(\int \kappa ds\right) = \frac{\bar{\kappa}}{\sqrt{\bar{\kappa}^2 + \bar{\tau}^2}}, \tag{20}$$

$$\sin\left(\int \kappa ds\right) = \frac{\bar{\tau}}{\sqrt{\bar{\kappa}^2 + \bar{\tau}^2}}, \tag{21}$$

respectively. Differentiating (21) with respect to $s$ we have

$$\kappa \cos\left(\int \kappa ds\right) = \frac{\bar{\tau}' \bar{\kappa}^2 - \bar{\tau} \bar{\kappa} \bar{\kappa}'}{(\bar{\kappa}^2 + \bar{\tau}^2)^{3/2}}. \tag{22}$$

From (14), (19) and (22), it follows

$$\kappa = \frac{\overline{\tau}'\overline{\kappa} - \overline{\tau}\,\overline{\kappa}'}{\overline{\kappa}^2 + \overline{\tau}^2},$$

or equivalently,

$$\kappa = \frac{\overline{\kappa}^2}{\overline{\kappa}^2 + \overline{\tau}^2}\left(\frac{\overline{\tau}}{\overline{\kappa}}\right)', \qquad (23)$$

which completes proof.

**Corollary 3.2.** *Let $\gamma$ with the curvature $\overline{\kappa}$ and the torsion $\overline{\tau}$ be an osculating-direction curve of $\alpha$. Then*

$$\frac{\kappa}{\tau} = \frac{\overline{\kappa}^2}{\left(\overline{\kappa}^2 + \overline{\tau}^2\right)^{3/2}}\left(\frac{\overline{\tau}}{\overline{\kappa}}\right)', \qquad (24)$$

*is satisfied, where $\kappa$ and $\tau$ are curvature and torsion of curve $\alpha$, respectively.*

### 4. Applications of osculating-direction curves

In this section, we give relationships between the osculating-direction curves and some special curves such as general helix, slant helix or rectifying curve in $E^3$.

*4.1. General helices and slant helices*

Considering Corollary 3.2, we have the following theorems which gives a way to construct the examples of slant helices by using general helices.

**Theorem 4.1.** *Let $\alpha: I \to \mathbb{R}^3$ be a Frenet curve in $E^3$ and $\gamma$ be an osculating-direction curve of $\alpha$. Then the followings are equivalent,*

  *i) A Frenet curve $\alpha$ is a general helix in $E^3$.*
  *ii) $\alpha$ is an osculating-donor curve of a slant helix.*
  *iii) An osculating-direction curve of $\alpha$ is a slant helix.*

**Theorem 4.2.** *Let $\alpha: I \to \mathbb{R}^3$ be a Frenet curve in $E^3$. An osculating-direction curve of $\alpha$ cannot be a general helix.*

**Example 4.1.** Let consider the general helix $\alpha(s) = \left(\cos\frac{s}{\sqrt{2}}, \sin\frac{s}{\sqrt{2}}, \frac{s}{\sqrt{2}}\right)$ in $E^3$ (Fig1). The Frenet vectors and curvatures of $\alpha$ are obtained as follows,

$$T(s) = \left(-\frac{1}{\sqrt{2}}\sin\frac{s}{\sqrt{2}},\ \frac{1}{\sqrt{2}}\cos\frac{s}{\sqrt{2}},\ \frac{1}{\sqrt{2}}\right),$$

$$N(s) = \left(-\cos\frac{s}{\sqrt{2}},\ \sin\frac{s}{\sqrt{2}},\ 0\right),$$

$$B(s) = \left(\frac{1}{\sqrt{2}}\sin\frac{s}{\sqrt{2}},\ -\frac{1}{\sqrt{2}}\cos\frac{s}{\sqrt{2}},\ \frac{1}{\sqrt{2}}\right),$$

$$\kappa = \tau = \frac{1}{2}.$$

Then we have $X(s) = (x_1(s), x_2(s), x_3(s))$ where

$$x_1(s) = -\frac{1}{\sqrt{2}}\sin\left(\frac{s}{2}+c\right)\sin\frac{s}{\sqrt{2}} - \cos\left(\frac{s}{2}+c\right)\cos\frac{s}{\sqrt{2}},$$

$$x_2(s) = \frac{1}{\sqrt{2}}\sin\left(\frac{s}{2}+c\right)\cos\frac{s}{\sqrt{2}} + \cos\left(\frac{s}{2}+c\right)\sin\frac{s}{\sqrt{2}},$$

$$x_3(s) = \frac{1}{\sqrt{2}}\sin\left(\frac{s}{2}+c\right).$$

and $c$ is integration constant. Now, we can construct a slant helix $\gamma$ which is also an osculating-direction curve of $\alpha$ (Fig. 2):

$$\gamma = \int_0^s \gamma'(s)ds = \int_0^s X(s)ds = (\gamma_1(s), \gamma_2(s), \gamma_3(s)),$$

where

$$\gamma_1(s) = \int_0^s \left[-\frac{1}{\sqrt{2}}\sin\left(\frac{s}{2}+c\right)\sin\frac{s}{\sqrt{2}} - \cos\left(\frac{s}{2}+c\right)\cos\frac{s}{\sqrt{2}}\right]ds,$$

$$\gamma_2(s) = \int_0^s \left[\frac{1}{\sqrt{2}}\sin\left(\frac{s}{2}+c\right)\cos\frac{s}{\sqrt{2}} + \cos\left(\frac{s}{2}+c\right)\sin\frac{s}{\sqrt{2}}\right]ds,$$

$$\gamma_3(s) = \int_0^s \frac{1}{\sqrt{2}}\sin\left(\frac{s}{2}+c\right)ds.$$

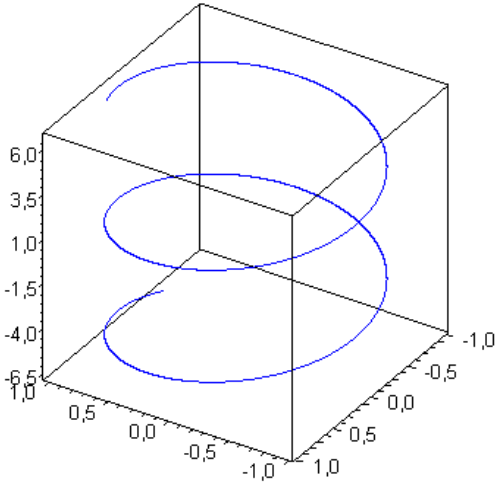
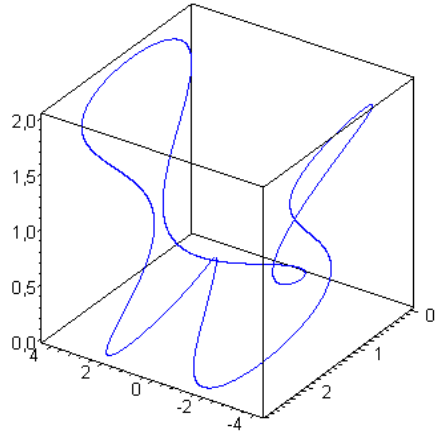

**Fig. 1.** General helix  **Fig. 2.** A slant helix constructed by a general helix.

**Example 4.2.** Let consider the general equation of spherical helix in $E^3$ defined by

$$\alpha(s) = \left(\cos s \cos(ws) + \frac{1}{w}\sin s \sin(ws),\ -\cos s \sin(ws) + \frac{1}{w}\sin s \cos(ws),\ \frac{\sin s}{cw}\right),$$

where $w = \frac{\sqrt{1+c^2}}{c}$, $c \in \mathbb{R}$ are constants. (See [14] (Fig3). The Frenet vectors and curvatures of $\alpha$ are given by the following equalities,

$$T(s) = \frac{c}{\sqrt{1+c^2(1-w^2)^2}} \left( (1-w^2)\sin(ws),\ (1-w^2)\cos(ws),\ \frac{1}{c} \right),$$

$$N(s) = \left( \cos(ws),\ -\sin(ws),\ 0 \right),$$

$$B(s) = \frac{1}{\sqrt{1+c^2(1-w^2)^2}} \left( \sin(ws),\ \cos(ws),\ -c(1-w^2) \right),$$

$$\kappa = \frac{c^2 w^2 (1-w^2)}{1+c^2(1-w^2)^2} \frac{1}{\cos s},$$

$$\tau = \frac{cw^2}{1+c^2(1-w^2)^2} \frac{1}{\cos s}.$$

Then we have $X(s) = (x_1(s),\ x_2(s),\ x_3(s))$ where

$$x_1(s) = \frac{c(1-w^2)}{\sqrt{1+c^2(1-w^2)^2}} \sin\left( A\ln(\sec s + \tan s) + c_1 \right) \sin(ws) + \cos\left( A\ln(\sec s + \tan s) + c_1 \right) \cos(ws),$$

$$x_2(s) = \frac{c(1-w^2)}{\sqrt{1+c^2(1-w^2)^2}} \sin\left( A\ln(\sec s + \tan s) + c_1 \right) \cos(ws) - \cos\left( A\ln(\sec s + \tan s) + c_1 \right) \sin(ws),$$

$$x_3(s) = \frac{1}{\sqrt{1+c^2(1-w^2)^2}} \sin\left( A\ln(\sec s + \tan s) + c_1 \right),$$

where $A = \dfrac{c^2 w^2 (1-w^2)}{1+c^2(1-w^2)^2}$ and $c_1$ is integration constant. Now, we can construct a slant helix $\gamma$ which is also an osculating-direction curve of $\alpha$ (Fig. 4):

$$\gamma = \int_0^s \gamma'(s)ds = \int_0^s X(s)ds = (\gamma_1(s),\ \gamma_2(s),\ \gamma_3(s)),$$

where

$$\gamma_1(s) = \int_0^s \left[ \frac{c(1-w^2)}{\sqrt{1+c^2(1-w^2)^2}} \sin\left( A\ln(\sec s + \tan s) + c_1 \right) \sin(ws) + \cos\left( A\ln(\sec s + \tan s) + c_1 \right) \cos(ws) \right] ds,$$

$$\gamma_2(s) = \int_0^s \left[ \frac{c(1-w^2)}{\sqrt{1+c^2(1-w^2)^2}} \sin\left( A\ln(\sec s + \tan s) + c_1 \right) \cos(ws) - \cos\left( A\ln(\sec s + \tan s) + c_1 \right) \sin(ws) \right] ds,$$

$$\gamma_3(s) = \int_0^s \frac{1}{\sqrt{1+c^2(1-w^2)^2}} \sin\left( A\ln(\sec s + \tan s) + c_1 \right) ds.$$

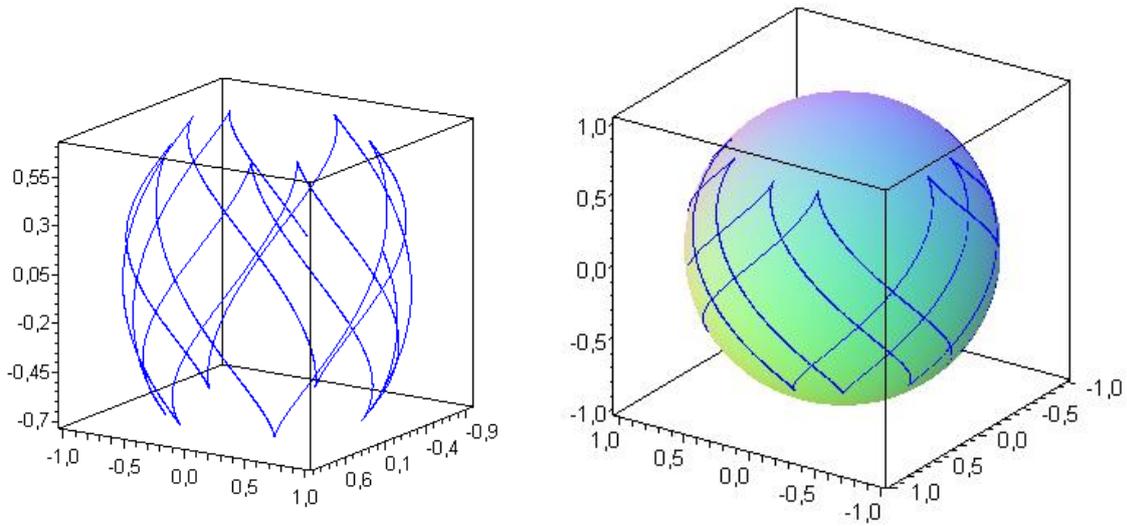

**Fig. 3.** Spherical helix

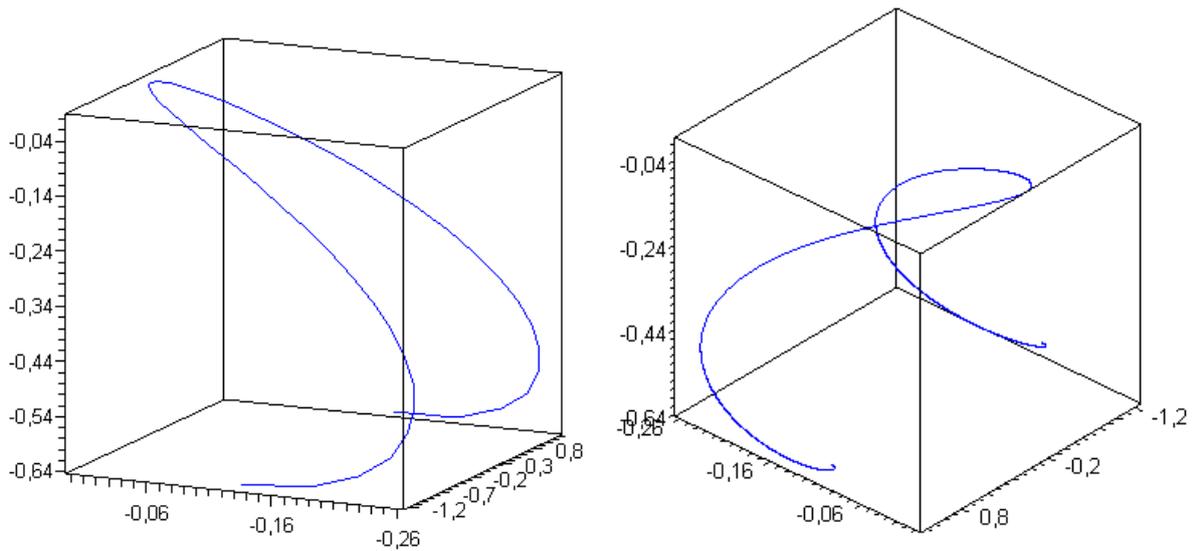

**Fig. 4.** A slant helix constructed by a spherical helix.

*4.2. Plane curves*

In this subsection, we give the relationship between osculating-direction curve and a plane curve. From Theorem 3.3 we can give the following theorem. Considering (11) and (23), the proof is clear, since $\alpha: I \to \mathbb{R}^3$ is a Frenet curve, i.e., $\kappa \neq 0$ in $E^3$.

**Theorem 4.3.** *Let $\alpha: I \to \mathbb{R}^3$ be a Frenet curve in $E^3$ and $\gamma$ be an osculating-direction curve of $\alpha$. Then the followings are equivalent,*

   *i) A Frenet curve $\alpha$ is not a plane curve in $E^3$.*
   *ii) An osculating-direction curve $\gamma$ of $\alpha$ is not a line or a plane curve in $E^3$.*
   *iii) $\alpha$ is not an osculating donor curve of a line or a plane curve in $E^3$*

*4.3. OD-Osculating Curves*

In this subsection we define osculating-direction ($OD$)-osculating curves in $E^3$ and give the relationships between osculating direction curves and $OD$-osculating curves.

A space curve whose position vector always lies in its rectifying plane is called rectifying curve [5]. Moreover, if the Frenet frame and curvatures of a space curve are given by $\{T, N, B\}$ and $\kappa, \tau$ respectively, then the vector $\tilde{D}(s) = \frac{\tau}{\kappa}(s)T(s) + B(s)$ is called modified Darboux vector of the curve [10].

Let now $\alpha$ be a Frenet curve with Frenet frame $\{T, N, B\}$ and $\gamma$ an osculating-direction curve of $\alpha$. The curve $\gamma$ is called osculating-direction osculating curve (or $OD$-osculating curve) of $\alpha$, if the position vector of $\gamma$ always lies on the osculating plane of its osculating-donor curve $\alpha$.

The definition of $OD$-osculating curve allows us to write the following equality,
$$\gamma(s) = m(s)T(s) + n(s)N(s) \tag{28}$$
where $m(s), n(s)$ are non-zero differentiable functions of $s$. Since $\gamma$ is osculating-direction curve of $\alpha$, from Corollary 3.4, we have
$$\begin{cases} T = \sin\left(\int \kappa ds\right)\bar{T} + \cos\left(\int \kappa ds\right)\bar{B}, \\ N = \cos\left(\int \kappa ds\right)\bar{T} - \sin\left(\int \kappa ds\right)\bar{B}. \end{cases} \tag{29}$$
Substituting (29) in (28) gives
$$\gamma(s) = \left[m\sin\left(\int \kappa ds\right) + n\cos\left(\int \kappa ds\right)\right]\bar{T} + \left[m\cos\left(\int \kappa ds\right) - n\sin\left(\int \kappa ds\right)\right]\bar{B}. \tag{30}$$
Writing
$$\begin{cases} \rho(s) = m\sin\left(\int \kappa ds\right) + n\cos\left(\int \kappa ds\right), \\ \sigma(s) = m\cos\left(\int \kappa ds\right) - n\sin\left(\int \kappa ds\right), \end{cases} \tag{31}$$
in (30) and differentiating the obtained equality we obtain
$$\bar{T} = \rho'\bar{T} + (\rho\bar{\kappa} - \sigma\bar{\tau})\bar{N} + \sigma'\bar{B}. \tag{32}$$
Then we have
$$\sigma = a = \text{constant}, \quad \rho = s + b = \frac{\bar{\tau}}{\bar{\kappa}}a, \tag{33}$$
where $a, b$ are non-zero integration constants. From (33) it follows that
$$\gamma(s) = a\left(\frac{\bar{\tau}}{\bar{\kappa}}\bar{T} + \bar{B}\right)(s) = a\tilde{\bar{D}}(s), \tag{34}$$
where $\tilde{\bar{D}}$ is the modified Darboux vector of $\gamma$.

Now we can give the followings which characterize $OD$-osculating curves.

**Theorem 4.4.** *Let $\alpha: I \to \mathbb{R}^3$ be a Frenet curve in $E^3$ and $\gamma$ be an osculating-direction curve of $\alpha$. If $\gamma$ is an $OD$-osculating curve in $E^3$, then we have the followings,*

  i) *$\gamma$ is a rectifying curve in $E^3$ whose curvatures satisfy $\frac{\bar{\tau}}{\bar{\kappa}} = \frac{s+b}{a}$ where $a, b$ are non-zero constants.*

*ii) The position vector and modified Darboux vector $\tilde{D}$ of $\gamma$ are linearly dependent.*

Theorem 4.4 gives a way to construct a rectifying curve by using osculating-donor curve as follows:

**Corollary 4.1.** *Let $\alpha: I \to \mathbb{R}^3$ be a Frenet curve in $E^3$ and $\gamma$ an OD-osculating curve of $\alpha$ in $E^3$. Then the position vector of $\gamma$ is obtained as follows,*

$$\gamma(s) = \left[(s+b)\sin\left(\int \kappa ds\right) + a\cos\left(\int \kappa ds\right)\right] T(s) + \left[(s+b)\cos\left(\int \kappa ds\right) - a\sin\left(\int \kappa ds\right)\right] N(s) \quad (35)$$

*where $a$, $b$ are non-zero integration constants.*
**Proof.** The proof is clear from (31), (33) and (28).

**Example 4.3.** Let consider the Frenet curve $\alpha: I \to \mathbb{R}^3$ with the parametrization

$$\alpha(s) = \left(\frac{\sqrt{2}}{3} s^{3/2}, \frac{\sqrt{2}}{3}(1-s)^{3/2}, \frac{\sqrt{2}}{2} s\right), \quad s \neq 0, \ s < 1.$$

Frenet vectors and curvatures of $\alpha$ are obtained as follows,

$$T(s) = \frac{\sqrt{2}}{2}\left(\sqrt{s}, -\sqrt{1-s}, 1\right),$$

$$N(s) = \left(\sqrt{1-s}, \sqrt{s}, 0\right),$$

$$B(s) = \frac{\sqrt{2}}{2}\left(-\sqrt{s}, \sqrt{1-s}, 1\right),$$

$$\kappa = \frac{\sqrt{2}}{4} \frac{1}{\sqrt{s(1-s)}}, \quad \tau = \frac{1}{2\sqrt{2s(1-s)}}.$$

By choosing $a = b = 1$, from (35) a rectifying curve $\gamma$ is given by the parametrization $\gamma(s) = (\gamma_1, \gamma_2, \gamma_3)$ where

$$\gamma_1(s) = \frac{\sqrt{2}}{2}\left[(s+1)\sin\left(\frac{\sqrt{2}}{4}\arcsin(2s-1) + d\right) + \cos\left(\frac{\sqrt{2}}{4}\arcsin(2s-1) + d\right)\right]\sqrt{s}$$

$$+ \left[(s+1)\cos\left(\frac{\sqrt{2}}{4}\arcsin(2s-1) + d\right) - \sin\left(\frac{\sqrt{2}}{4}\arcsin(2s-1) + d\right)\right]\sqrt{1-s},$$

$$\gamma_2(s) = -\frac{\sqrt{2}}{2}\left[(s+1)\sin\left(\frac{\sqrt{2}}{4}\arcsin(2s-1) + d\right) + \cos\left(\frac{\sqrt{2}}{4}\arcsin(2s-1) + d\right)\right]\sqrt{1-s}$$

$$+ \left[(s+1)\cos\left(\frac{\sqrt{2}}{4}\arcsin(2s-1) + d\right) - \sin\left(\frac{\sqrt{2}}{4}\arcsin(2s-1) + d\right)\right]\sqrt{s},$$

$$\gamma_3(s) = \frac{\sqrt{2}}{2}\left[(s+1)\sin\left(\frac{\sqrt{2}}{4}\arcsin(2s-1) + d\right) + \cos\left(\frac{\sqrt{2}}{4}\arcsin(2s-1) + d\right)\right],$$

and $d$ is integration constant (Fig. 6).

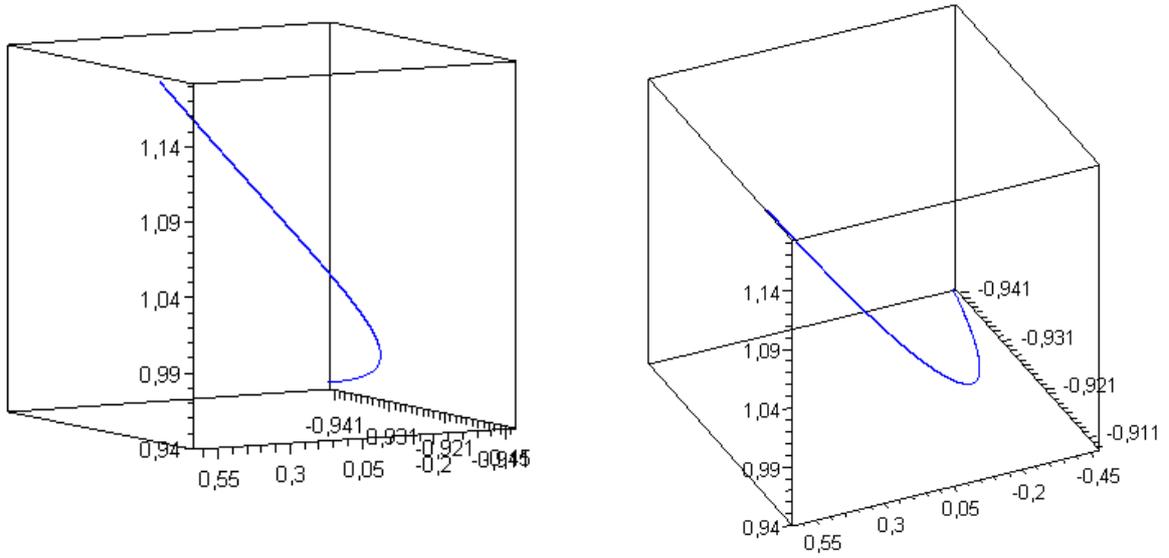

**Fig. 6.** $OD$-osculating curve

**Example 4.4.** Let consider the Frenet curve $\alpha : I \to \mathbb{R}^3$ with the parametrization
$$\alpha(s) = \left(12\cos\frac{s}{13},\ 12\sin\frac{s}{13},\ \frac{5s}{13}\right).$$
Frenet vectors and curvatures of $\alpha$ are obtained as follows,
$$T(s) = \frac{\alpha'}{\|\alpha'\|} = \left(\frac{-12}{13}\sin\frac{s}{13},\ \frac{12}{13}\cos\frac{s}{13},\ \frac{5}{13}\right),$$
$$N(s) = \left(-\cos\frac{s}{13},\ -\sin\frac{s}{13},\ 0\right),$$
$$B(s) = \left(\frac{5}{13}\sin\frac{s}{13},\ -\frac{5}{13}\cos\frac{s}{13},\ \frac{12}{13}\right),$$
$$\kappa = \frac{12}{169},\quad \tau = \frac{5}{169}.$$
By choosing $a = b = 1$, from (35) a rectifying curve $\gamma$ is given by the parametrization $\gamma(s) = (\gamma_1, \gamma_2, \gamma_3)$ where
$$\gamma_1(s) = -\frac{12}{13}\left[(s+1)\sin\left(\frac{12s}{169}+e\right)+\cos\left(\frac{12s}{169}+e\right)\right]\sin\frac{s}{13}$$
$$-\left[(s+1)\cos\left(\frac{12s}{169}+e\right)-\sin\left(\frac{12s}{169}+e\right)\right]\cos\frac{s}{13},$$
$$\gamma_2(s) = \frac{12}{13}\left[(s+1)\sin\left(\frac{12s}{169}+e\right)+\cos\left(\frac{12s}{169}+e\right)\right]\cos\frac{s}{13}$$
$$-\left[(s+1)\cos\left(\frac{12s}{169}+e\right)-\sin\left(\frac{12s}{169}+e\right)\right]\sin\frac{s}{13},$$
$$\gamma_3(s) = \frac{5}{13}\left[(s+1)\sin\left(\frac{12s}{169}+e\right)+\cos\left(\frac{12s}{169}+e\right)\right],$$

and *e* is integration constant (Fig. 7).

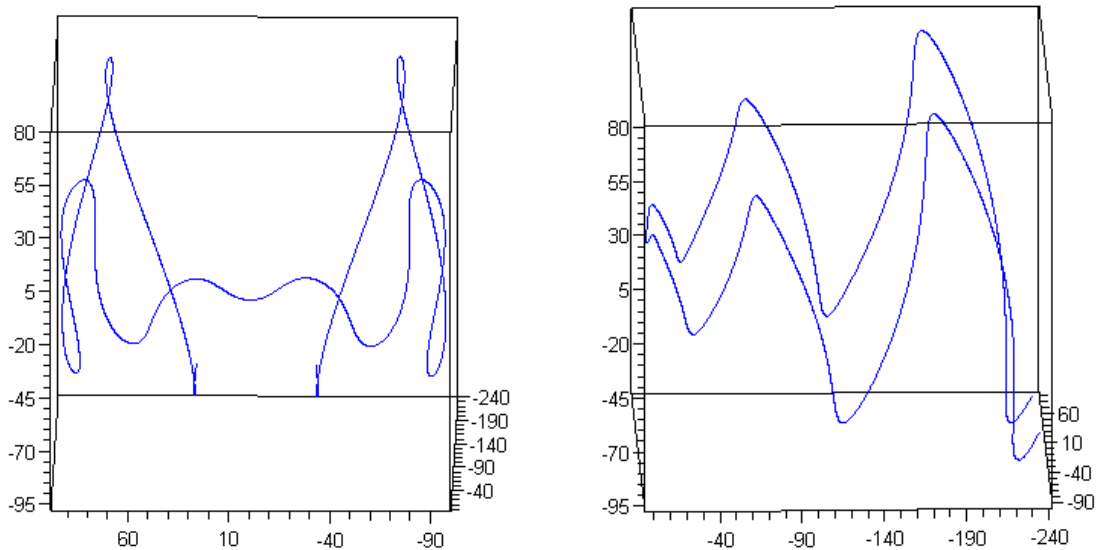

**Fig. 7.** *OD*-osculating curve

[15] Otsuki, T., Differential Geometry, Asakura Publishing Co. Ltd., Tokyo, (in Japanese), (1961).

[16] Puig-Pey, J. Gálvez, A., Iglesias, A., Helical curves on surfaces for computer-aided geometric design and manufacturing, in: Computational Science and its Applications-ICCSA 2004. Part II, 771-778, in: Lecture Notes in Comput. Sci., vol. 3044, Springer, Berlin, (2004).

[17] Struik, D.J., Lectures on Classical Differential Geometry, $2^{nd}$ ed. Addison Wesley, Dover, (1988).

[18] Wang, F., Liu, H., Mannheim partner curves in 3-Euclidean space, Mathematics in Practice and Theory, 37(1) (2007) 141-143.

[19] Whittemore, J.K., Bertrand curves and helices. Duke Math. J. 6(1) (1940) 235-245.

[20] Wong, Y.C., On the generalized helices of Hayden and Syptak in an $N$-space, Proc. Cambridge Philos. Soc. 37 (1941) 229-243.